\documentclass[10pt]{article}

\usepackage{hyperref}
\usepackage[utf8]{inputenc}
\usepackage[OT2,T2A,T1]{fontenc}
\usepackage[english,french]{babel}
\usepackage[osf]{mathpazo}

\usepackage[pdftex]{graphicx}

\usepackage{amsmath,amsfonts,amsthm,amssymb}

\usepackage{color}
\definecolor{FOB}{rgb}{0.,0.,0.7}
\definecolor{FOR}{rgb}{0.65,0.,0.}
\definecolor{FOfond}{rgb}{0.98,0.98,0.91}

\pagecolor{FOfond}

\usepackage{sectsty}
\allsectionsfont{\color{FOR}}

\def\rd{{\mathrm d}}

\def\cA{{\mathcal A}}
\def\cB{{\mathcal B}}

\def\cF{{\mathcal F}}

\def\cK{{\mathcal K}}

\def\cM{{\mathcal M}}

\def\cP{{\mathcal P}}
\def\cQ{{\mathcal Q}}
\def\cR{{\mathcal R}}

\def\ord{{\rm ord}}

\def\K{{\bf K}}
\def\N{{\bf N}}

\def\Q{{\bf Q}}

\def\bull{\vrule height .9ex width .8ex depth -.1ex }
\def\rd{{\rm d}}

\def\Spol{\hbox{\rm S-Pol}}

\def\eqref #1{(\ref{#1})}

\newcounter{thenum}
\def\texttheo{\relax}
\newenvironment{theorem}{\medbreak\refstepcounter{thenum}
\noindent\textsc{Theorem} %
\thethenum. \texttheo ---  \it  }{\rm }
\newenvironment{e-proposition}{\medbreak\refstepcounter{thenum}
\noindent\textsc{Proposition} \thethenum. ---  \it  }{\rm }

\newenvironment{e-definition}{\medbreak\refstepcounter{thenum}
\noindent\textsc{Definition} \thethenum. ---  \it  }{\rm }
\newenvironment{lemma}{\medbreak\refstepcounter{thenum}\noindent{\it Lemma} %
\thethenum. ---  \it  }{\rm }

\newenvironment{e-rem}{\medbreak\refstepcounter{thenum}{}%
 \thethenum) }{}

\newenvironment{e-ex}{\medbreak\refstepcounter{thenum}{}%
 \thethenum) }{}

\newenvironment{theoreme}{\medbreak\refstepcounter{thenum}\noindent\textsc{Théorème} %
\thethenum. \texttheo---  \it  }{\rm }

\newenvironment{definition}{\medbreak\refstepcounter{thenum}
\noindent\textsc{Définition} \thethenum. ---  \it  }{\rm }
\newenvironment{lemme}{\medbreak\refstepcounter{thenum}\noindent{\it Lemme} %
\thethenum. ---  \it  }{\rm }

\newenvironment{exemple}{\medbreak\refstepcounter{thenum}\noindent{\sl
    Exemple} %
\thethenum. --- }{}

\newenvironment{preuve}{\smallbreak\noindent{\sc Preuve.} ---
  \rm}{\quad\bull\smallskip\rm}

\begin{document}

\thispagestyle{empty}

\hbox to\hsize{\hss\vbox{\hsize=15cm
\begin{center}{{\huge\parindent=0pt\selectlanguage{french} 
      \color{FOR}
      Une généralization du critère

      de Boulier -- Buchberger

      pour le calcul des ensembles caractéristiques

      d'idéaux différentiels

      }}
\vskip 0.5cm

{{\Large\parindent=0pt\selectlanguage{english}
      \color{FOR}
    A generalization of the Boulier -- Buchberger criterion

    for the computation of characteristic sets

    of differential ideals
}}
\end{center}}\hss}
\vfill

\hbox to \hsize{\vbox{\hsize = 5.8cm\parindent=0pt
    {\large Amir \textsc{Hashemi}}
\bigskip

Department of Mathematical Sciences,

Isfahan University of Technology,

Isfahan, 84156-83111,

Iran
\smallskip

{\scriptsize Amir.Hashemi@iut.ac.ir}

  }

  \hss\vbox{\hsize = 5.5cm \parindent=0pt

    {\large François \textsc{Ollivier}}
\bigskip

CNRS: LIX, UMR 7161 

École polytechnique 

91128 Palaiseau Cedex

France

\smallskip

{\scriptsize francois.ollivier@lix.polytechnique.fr}

}}
\vfill

\begin{center}{\parindent=0pt\large Mars 2022

}
\end{center}
\vfill

{\small
\hbox to \hsize{\vbox{\hsize=
    5.8cm\selectlanguage{french}
    \spaceskip=.3333em plus.3333em minus.1111em
\parindent=0pt \selectlanguage{french}\small
\noindent \textbf{Résumé}  Nous généralisons l'analogue du premier
critère de Buchberger, dû à Boulier \textit{et al.}, pour détecter les
réductions inutiles de S-poly\-nômes, lors des calculs d'ensembles
caractéristiques d'idéaux différentiels. La version primitive suppose
des polynômes linéaires; le résultat est ici étendu à un produit de
polynômes différentiels linéaires, appliqués à un même polynôme
différentiel, arbitraire.

}
  \hss

  \vbox{\hsize= 5.5cm\selectlanguage{english}
    \spaceskip=.3333em plus.3333em minus.1111em
\parindent=0pt\small
\noindent \textbf{Abstract.}  We generalize the analog of Buchberger's first criterion, stated by
Bou\-lier \textit{et al.}, for detecting useless S-polynomials reductions in the
computation of characteristic sets of differential ideals. The original
version assumes linear polynomials; this result is here extended to a
product of linear differential polynomials depending of the same
arbitrary differential polynomial.

}}}
\vskip 0.5cm

\noindent Classification AMS: 12H05

\vskip 0.5cm

\hbox to\hsize{\vbox{\noindent\hsize=5.8cm {\small {\bf Mots-clés:}
      Algèbre différentielle, Ensembles caractéristique,
  Algorithme Rosenfeld--\hfill\break Gröbner, Critère de Boulier--Buchberger
}}\hss

\vbox{\noindent\hsize=5.5cm
    {\small\selectlanguage{english} {\bf Keywords:} Differential
      algebra, Characteristic sets, Rosenfeld--Gröbner algorithm,
      Boulier--Buchberger criterion}}}
\medskip

\noindent{\tiny Version d'auteur de l'article:
  \textsc{Hashemi} (Amir) et \textsc{Ollivier} (François), \og Une
  généralization du critère de Boulier -- Buchberger pour le calcul
  des ensembles caractéristiques d'idéaux différentiels \fg,
  \textit{Comptes Rendus. Mathématique}, Tome~360, 255--264, 2022.
 \hfill
 DOI:~\href{https://doi.org/10.5802/crmath.295}{10.5802/crmath.295}\break

}

\vfill\eject

\selectlanguage{english}
 \section*{Abridged English version}
 {\def\notasname{Notation}%

 Boulier \textit{et al.}~\cite[Prop.~4 p.~92]{Boulier2009} gave a
 differential analog of Buchberger's first
 criterion~\cite{Buchberger1979}, stating that if leading monomials of
 two polynomials have no common factor, then their S-polynomial will
 be reduced to zero by these two polynomials. The original
 differential version requires linear polynomials, a result that was
 extended to products of linear factor in Hashemi and
 Touraji~\cite{Hashemi2014}, without a complete proof. In this paper
 we prove a more general statement, where the linear factors are no more
 assumed to depend on the same variable, but only on the same
 differential polynomial $z$.

  We consider the differential algebra $\cR:=\cF\{X\}$, where
  $X:=\{x_{1}, \ldots, x_{n}\}$ and $\cF$ is a differential field of
  characteristic~$0$, with a set of mutually commuting derivations
  $\Delta=\{\delta_{1}, \ldots, \delta_{m}\}$. We denote by $\Theta$
  the commutative monoïd generated by $\Delta$ and $\Upsilon:=\Theta
  X$ is the set of derivatives. The leading derivative of a
  differential polynomial $P\in\cR$ is written $\upsilon_{P}$ and $[1,
    k]:=\{1, \ldots, k\}\subset\N$. Let $D\subset\Delta$, we denote by
  $\cF_{D}$ the differential field restricted to derivations
  $\delta\in D$ and by $\cR_{D}:=\cF_{D}\{x\}$ the corresponding
  differential ring. In the same way, $\Theta_{D}$ is the commutative
  monoïde generated by $D$ and $\K_{D}\subset\cF$ the field of
  constants for derivations in $D$.

  We refer to Kolchin~\cite{Kolchin1973} and Boulier \textit{et
    al.}~\cite{Boulier2009} 
  for basic definitions and notations of differential algebra. We
  denote by $I_{P_{i}}$ or $I_{i}$ and $S_{P_{i}}$ or $S_{i}$ the
  initial and separant of a polynomial $P_{i}$, by $I_{\cA}$
  (resp.\ $S_{\cA}$) the product of initials (resp.\ separants) of
  polynomials in a set $\cA$ and $H_{\cA}:=I_{\cA}S_{\cA}$. The
  notations $(\Sigma)_{R}$ or $[\Sigma]_{R}$ denote the algebraic or
  differential ideals generated by $\Sigma$ in the ring $R$

 We have then this first generalization of Boulier's theorem, using
 linear operators  $\theta_{i}-L_{i}$, applied to a differential
 polynomial $z$, that is only required not to belong to the ground
 field.

 \begin{theorem}\label{th:eBoulier} Let $\Delta_{i}\subset\Delta$,
  $i=1, 2, 3$, with $\Delta_{i}\cap\Delta_{j}=\emptyset$ for $i\neq j$
   and $P_{i}\in\cR$, $i=1, 2$. Let $z\in\cR\setminus\cF$ and let
   $\prec$ be an admissible order on derivatives, such that
   $\upsilon_{P_{i}}=\theta_{i}\upsilon_{z}$, with
   $1\neq\theta_{i}\in\Theta_{\Delta_{i}}$. We further assume that
   $P_{i}=(\theta_{i}-L_{i})z$, where $L_{1}L_{2}=L_{2}L_{1}$ and
   $L_{i}\in\K_{\Delta_{j}}[\Delta_{i}\cup \Delta_{3}]$, where $\{i,
   j\}=\{1, 2\}$. Under those hypotheses, using $\prec$ to define
   reduction, we have the following propositions.
   
i)  The polynomial
$\Spol(P_{1},P_{2})=\theta_{2}P_{1}-\theta_{1}P_{2}$
is reduced to $0$ by $\{P_{1},P_{2}\}$, using differential standard
bases reduction (\textit{cf.}~\cite{Carra1989,Ollivier1990}).

ii) It is reduced to $0$ using characteristic sets pseudo-reduction.

iii) If $z$ is linear, the set $\{P_{1},P_{2}\}$ is a differential
standard basis of the prime differential ideal $[P_{1},P_{2}]$.

iv) The set $\{P_{1},P_{2}\}$ is a characteristic set of the prime
differential ideal
$[P_{1},P_{2}]:S_{z}^{\infty}=[A_{1},A_{2}]:S_{A}^{\infty}$, where
$A_{i}$ is the factor of $P_{i}$, the main derivative of which is
$\theta_{i}\upsilon_{z}$. 
\end{theorem}

The next lemma extends the assertion~ii) and iv) above to square-free
products.
 
 \begin{lemma}\label{lem:eprod}
  Let $P_{1, k}$, for $1\le k\le r_{1}$, and $P_{2, k}$, for $1\le
  k\le r_{2}$, be two finite families of irreducible polynomials, of
  degree $1$ in their leading derivatives and such that the $P_{1,k}$
  (resp.~$P_{2,k}$) have the same leading derivative $\upsilon_{1}$
  (resp.~$\upsilon_{2}$) and do not depend on $\upsilon_{2}$
  (resp.~$\upsilon_{1}$). We assume that the prime differential ideals
  defined by those polynomials (according to th.~\ref{th:eBoulier} iv)
  are mutually different. Let $Q_{i}:=\prod_{k=1}^{r_{1}} P_{1, k}$,
  for $i=1,2$, and $Q:=\{Q_{1},Q_{2}\}$.
   \smallskip

  i) Under those hypotheses, the following propositions are equivalent:

\noindent a) for all $(j_{0}, k_{0})\in [1,r_{1}]\times[1,r_{2}]$,
$\Spol(P_{1, j_{0}},P_{2,  k_{0}})$ is peudo-reduced to $0$ by
$\{P_{1, j_{0}}, P_{2,  k_{0}}\}$; 

\noindent b) $\Spol(Q_{1},Q_{2})$
  is pseudo-reduced to $0$ by $Q$.
   \smallskip

  ii) If a) or b) stands, then $Q$ is a characteristic set and a
  characteristic representation of
  $[Q]:H_{Q}^{\infty}=
  \bigcap_{j=1}^{r_{1}}\bigcap_{k=1}^{r_{2}}[P_{1,j},P_{2,k}]:(I_{1,j}I_{2,k})^{\infty}$.
\end{lemma}

When powers may appear in the products, the following technical lemma
is needed. For short, we use the notation $\partial_{y}:=\partial/\partial y$.

\begin{lemma}\label{lem:ecarre} Let $\cP\subset\cR$ be a prime ideal.
  We denote by $\Omega_{\cR/\cF}$ the module of Kähler differentials
  and by $\rd$ the canonical derivation from $\cR$ to
  $\Omega_{\cR/\cF}$.

i) For all $r\in\N$ and $Q\in\cP^{r}$, we have: $\rd
Q\in\cP^{r-1}\Omega_{\cR/\cF}$.

ii) Let $\cA=\{A_{1}, \ldots, A_{s}\}$ be a characteristic set of a
prime ideal $\cP$, for an ordering $\prec$ and $d_{i}\in\N$, for $1\le
i\le s$,  arbitrary integers.

\noindent a) For all $A\in\cR$, all $d\in\N$ and all
$(\tau,\theta)\in\Theta^{2}$, we have $\partial_{\tau\!\!A}\, \theta\!\!A^{d}\in
          [A^{d-1}]$.

\noindent b) If $Q\in\cQ:=[A_{i}^{d_{i}}|1\le i\le
  s]:H_{\cA}^{\infty}$, we have: $\rd Q=\sum_{i=1}^{s}\sum_{k\in\N}c_{i,
  k}\rd A_{i}^{(k)}$, with $c_{i, k}\in
          [A_{\iota}^{\bar d_{\iota}}|1\le\iota\le
            s]:H_{\cA}^{\infty}$, where $\bar d_{\iota}:=d_{i}$ if
          $\iota\neq i$ and $\bar d_{\iota}:=d_{i}-1$ if
          $\iota=i$.

\noindent c) With the same notations as in b), we further denote by
$\Xi$ the set of derivatives appearing in all the differential
polynomials $A_{i}$, for $1\le i\le s$ and different from the leading
derivatives $\upsilon_{A_{j}}$, for all $1\le j\le s$. If the $A_{i}$
are all linear in their leading derivatives $\upsilon_{A_{i}}$ and if
$I_{i}\in\cF[\Xi]$, we have the equality
$\cQ\cap\cR_{0}=(A_{i}^{d_{i}}|1\le i\le s):I_{\cA}^{\infty}$, where
$\cR_{0}:=\cF[\Xi,\upsilon_{A_{1}}, \ldots, \upsilon_{A_{s}}]$. In
this case, $\{A_{i}^{d_{i}}|1\le i\le s\}$ is a characteristic set of
$[A_{i}^{d_{i}}|1\le i\le s]:H_{\cA}^{\infty}$.
\end{lemma}

Using these two lemmas, we can prove the main theorem.

 \begin{theorem}\label{th:eprincipal} Let
  $\bar Q_{i}=\prod_{k=0}^{r_{i}} P_{i, k}^{d_{i, k}}$ and
   $Q_{i}=\prod_{k=1}^{r_{i}} P_{i, k}^{d_{i, k}}$, $i=1, 2$, where
   for any couple $(P_{1, k},P_{2, k'})$, $(k, k')\in[1,
     r_{1}]\times[1, r_{2}]$, the hypotheses of th.~\ref{th:eBoulier}
   are satisfied for the same $\theta_{i}z$, and $\upsilon_{P_{i,
       0}}\prec \theta_{i}\upsilon_{z}$, $i=1, 2$. We further assume
   that the $P_{i, k}$ are mutually different.  Under these
   hypotheses, we have the following assertions.
  
  i) If $d_{i,  k}=1$, $i=1, 2$ and $1\le k\le
  r_{i}$, then
$\bar Q:=\{\bar Q_{1}, \bar Q_{2}\}$ is a characteristic set of
$[\bar Q]:S_{\bar Q}^{\infty}=[Q]:S_{Q}^{\infty}$.

ii) For any $d_{i, k}$, $ T:=\Spol(\bar Q_{1}, \bar Q_{2})=S_{\bar
  Q_{2}}\theta_{2}\bar Q_{1} - S_{\bar Q_{1}}\theta_{1}\bar Q_{2} $ is
pseudo-reduced to $0$ by $Q$.
\end{theorem}

}
 
\selectlanguage{french}

\section*{Introduction}

Boulier \textit{et al.} \cite[prop.~4 p.~92]{Boulier2009} ont donné un
analogue différentiel du premier critère de
Buchberger~\cite{Buchberger1979}, qui affirme que si les monômes de
tête de deux polynômes sont sans facteur commun, alors le S-polynôme
correspondant sera réduit à zéro.  Une généralisation de ce résultat
au cas de produits de facteurs linéaires a été énoncée par Hashemi et
Touraji~\cite{Hashemi2014}. Nous donnons ici une preuve complète dans
le cas où ces facteurs ne dépendent pas uniquement d'une même
variable, mais d'un même polynôme arbitraire $z$.

 Dans la suite, on considère un anneau de polynômes différentiels
 $\cR:=\cF\{X\}$, où $X:=\{x_{1}, \ldots, x_{n}\}$ et $\cF$ est un
 corps différentiel de caractéristique $0$, muni d'un ensemble de
 dérivations $\Delta=\{\delta_{1}, \ldots, \delta_{m}\}$ commutant
 entre elles. On note $\Theta$ le monoïde commutatif engendré par
 $\Delta$ et $\Upsilon:=\Theta X$ l'ensemble des dérivées. La dérivée
 de tête d'un polynôme $P\in\cR$ est notée $\upsilon_{P}$, tandis que
 $S_{P_{i}}$ ou $S_{i}$ et $I_{P_{i}}$ ou $I_{i}$ désignent
 respectivement le séparant et l'initial de $P_{i}$, avec
 $S_{\cA}:=\prod_{A\in\cA}S_{A}$, $I_{\cA}:=\prod_{a\in\cA}I_{A}$,
 $H_{\cA}:=S_{\cA}I_{\cA}$, et $[1, k]\subset\N$ est l'ensemble $\{1,
 \ldots, k\}$.  Les notations $(\Sigma)_{R}$ ou $[\Sigma]_{R}$
 désignent les idéaux algébrique et différentiel engendrés par
 $\Sigma$ dans l'anneau $R$.

Nous renvoyons à Kolchin~\cite{Kolchin1973} et Boulier \textit{et al.}
\cite{Boulier2009} pour les définitions de base en algèbre
différentielle. Nous aurons aussi besoin de la définition suivante.

\begin{definition} Soit $D\subset\Delta$, on note $\cF_{D}$ le corps
  différentiel resteint aux dérivations $\delta\in D$ et
  $\cR_{D}=\cF_{D}\{x\}$ l'anneau différentiel correspondant. On
  définit de même $\Theta_{D}$, le monoïde commutatif engendré par
  $D$. On désigne par $\K_{D}\subset\cF$ le corps de constantes pour
  les dérivations de $D$.
\end{definition}

\section{Une première généralisation}

On peut énoncer le théorème suivant, qui est une généralisation facile
de l'analogue du premier citère de Buchberger de Boulier
\textit{et.~al}~\cite[prop.~4 p.~92]{Boulier2009}, en passant du cas
où les $P_{i}$ sont linéaires à celui où ils sont obtenus en
appliquant un opérateur linéaire $\theta_{i}-L_{i}$ à un même polynôme
différentiel $z$, que l'on suppose uniquement ne pas appartenir au
corps de base, et qui peut donc être non linéaire. Par ailleurs, on
n'a pas besoin de supposer que les opérateurs $L_{i}$ sont à
coefficients constants, mais seulement les hypothèses minimales pour
obtenir le résultat, incluant la commutation de
$L_{1}$ et $L_{2}$, toujours acquise si $\Delta_{3}=\emptyset$.

\begin{theoreme}\label{th:Boulier} Soient $\Delta_{i}\subset\Delta$,
  $i=1, 2, 3$, tous deux à deux disjoints.
  On considère deux polynômes
  $P_{i}\in\cR$, $i=1, 2$. Soit $z\in\cR\setminus\cF$, et $\prec$ un
  ordre sur les dérivées tel que
  $\upsilon_{P_{i}}=\theta_{i}\upsilon_{z}$, avec
  $1\neq\theta_{i}\in\Theta_{\Delta_{i}}$, $1=1,2$. On suppose en outre que
  $P_{i}=(\theta_{i}-L_{i})z$, où $L_{1}L_{2}=L_{2}L_{1}$ et
  $L_{i}\in\K_{\Delta_{j}}[\Delta_{i}\cup \Delta_{3}]$,
  pour $\{i, j\}=\{1, 2\}$.

i)  Sous ces hypothèses, pour tout ordre sur les monômes de $\cR$
  compatible avec $\prec$, le S-polynôme
$$
\Spol(P_{1},P_{2})=\theta_{2}P_{1}-\theta_{1}P_{2}
$$
est réduit à $0$ par $\{P_{1},P_{2}\}$ au sens des bases standards
différentielles (cf.~\cite{Carra1989,Ollivier1990}).

ii) Il est réduit à $0$ au sens des ensembles caractéristiques.

iii) Si $z$ est linéaire, l'ensemble $\{P_{1},P_{2}\}$ est une base
standard de l'idéal différentiel premier $[P_{1},P_{2}]$.

iv) L'ensemble $\{P_{1},P_{2}\}$ est un ensemble caractéristique de
l'idéal différentiel premier
$[P_{1},P_{2}]:S_{z}^{\infty}=[A_{1},A_{2}]:S_{A}^{\infty}$, où
$A_{i}$ désigne le facteur irréductible de $P_{i}$ dont la dérivée
dominante est $\theta_{i}\upsilon_{z}$. 
\end{theoreme}
\begin{preuve} i) Ce polynôme est égal à
  $(\theta_{2}L_{1}-\theta_{1}L_{2})z=(L_{1}\theta_{2}-L_{2}\theta_{1})z$,
  car, avec nos hypothèses, $\theta_{1}$ et $L_{2}$, $\theta_{2}$ et
  $L_{1}$ commutent. La réduction du monôme de tête de $L_{1}\theta_{2}$ par
  $P_{2}$ remplace $\theta_{2}z$ par $L_{2}z$ et inversement, en
  permutant les indices. On n'a donc pas besoin de multiplier par
  $S_{z}$; c'est une réduction au sens des bases standards. On obtient donc
  $(L_{1}L_{2}-L_{2}L_{1})z=0$, puisque $L_{1}$ et $L_{2}$ commutent.

  ii) Est une conséquence directe de i).

  iii) Aussi,  puisque le seul S-polynôme
  possible dans le cas linéaire est réduit à $0$
  (\textit{Cf.}~\cite[th.~5]{Ollivier1990}). La linéarité entraîne la
  primalité de l'idéal.

  iv) Aussi, en utilisant le lemme de Rosenfeld
  \cite[th.~3]{Boulier2009} et en remarquant que
  $I_{P_{1}}=S_{P_{1}}=I_{P_{2}}=S_{P_{2}}=S_{z}$.  On peut alors
  évidement remplacer $P_{i}$ par $A_{i}$ puisque les deux diffèrent
  par un facteur qui ne contient pas la dérivée dominante, et divise
  donc le séparant. Comme les
  $P_{i}$ sont de degré $1$ en leurs dérivées dominantes, les $A_{i}$
  sont absolument irréductibles, ce qui entraîne la primalité de
  $(A):S_{A}^{\infty}$ et de l'idéal différentiel $[A]:S_{A}^{\infty}$.
\end{preuve}

\begin{exemple}\label{ex::1} On considère le corps
  $\cF:=\Q(y_{1},y_{2},y_{3})$, 
  muni des dérivations $\delta_{i}:=\partial/\partial y_{i}$, pour
  $i=1,2,3$. Soient $\Delta_{1}:=\{\delta_{1}\}$,
  $\Delta_{2}:=\{\delta_{2}\}$, $\Delta_{3}:=\{\delta_{3}\}$,
  $\theta_{1}=\delta_{1}^{3}$, $\theta_{2}=\delta_{2}^{3}$,
  $L_{1}=(\delta_{1}-y_{1})(\delta_{3}-y_{3})$ et
  $L_{2}=(\delta_{2}-y_{2})(\delta_{3}-y_{3})$. On se place sur
  $\cR:=\cF\{x\}$, muni de l'ordre sur les dérivées qui utilise
  d'abord l'ordre total de dérivation, puis l'ordre lexicographique
  avec $\delta_{1}>\delta_{2}>\delta_{3}$.

  En posant $z=x$, on a
  $P_{1}=\delta_{1}^{3}x+\delta_{1}\delta_{3}x
  -y_{3}\delta_{1}x-y_{1}\delta_{3}x+y_{1}y_{3}x$ et 
  $P_{2}=\delta_{2}^{3}x+\delta_{2}\delta_{3}x
  -y_{3}\delta_{2}x-y_{2}\delta_{3}x+y_{2}y_{3}x$. On a alors 
$$
  \delta_{2}^{3}P_{1}-\delta_{1}^{3}P_{2}=
  (\delta_{1}-y_{1})(\delta_{3}-y_{3})\delta_{2}^{3}x -
  (\delta_{2}-y_{2})(\delta_{3}-y_{3})\delta_{1}^{3}x,
$$
qui est réduit par $\{P_{1},P_{2}\}$ à
$$
(\delta_{1}-y_{1})(\delta_{3}-y_{3})(\delta_{2}-y_{2})(\delta_{3}-y_{3})x -
  (\delta_{2}-y_{2})(\delta_{3}-y_{3})(\delta_{1}-y_{1})(\delta_{3}-y_{3})x\\
=0.
$$
Dans ce cas, $\{P_{1},P_{2}\}$ est un ensemble caractéristique et
aussi une base standard de l'idéal qu'il engendre.

Le calcul est semblable en prenant, \textit{e.g.} $z=x^{3}$, mais il
devient indispensable de multiplier par le séparant $3x^{2}$ de
$x^{3}$, qui est aussi le séparant et l'initial de $\delta_{i}x^{3}$
et $P_{i}$, pour $i=1,2$, afin d'effectuer la réduction.
\end{exemple}

\section{Le cas des produits}

On peut énoncer le lemme général suivant, qui est essentiellement
suffisant pour établir notre théorème principal~\ref{th:principal}
dans le cas de produits sans carrés, c'est-à-dire le i) de ce théorème.

\begin{lemme}\label{lem:prod}
  Soient $P_{1, k}$, pour tout $1\le k\le r_{1}$ et $P_{2, k}$, pour
  tout $1\le k\le r_{2}$ deux familles finies de polynômes
  irréductibles et de degré $1$ en leurs dérivées de tête. Les $P_{1,
    k}$ (resp.~$P_{2, k}$) ont la même dérivée de tête $\upsilon_{1}$
  (resp.~$\upsilon_{2}$) et ne dépendent pas de $\upsilon_{2}$
  (resp. $\upsilon_{1}$). On suppose que les idéaux différentiels
  premiers dont ces polynômes sont les ensembles caractéristiques (par
  le th.~\ref{th:Boulier} iv) sont sont tous différents, ce qui
  revient à dire que leurs facteurs irréductibles contenant la dérivée
  de tête sont différents.  On note $Q_{i}:=\prod_{k=1}^{r_{1}} P_{1,
    k}$, pour $i=1,2$ et $Q:=\{Q_{1},Q_{2}\}$.

i) Sous ces hypothèses, les
  propositions suivantes sont équivalentes:

\noindent a) pour tout $(j_{0}, k_{0})\in [1,r_{1}]\times[1,r_{2}]$,
$\Spol(P_{1, j_{0}},P_{2, k_{0}})$ est pseudo-réduit à $0$ par
$\{P_{1, j_{0}}, P_{2, k_{0}}\}$; 

\noindent b) $\Spol(Q_{1}, Q_{2})$
  est pseudo-réduit à $0$ par $\{Q_{1}, Q_{2}\}$.

  ii) Si a) ou b) est vérifiée, alors $Q$ est un
  ensemble caractéristique et une
  représentation caractéristique de
\begin{equation}\label{eq::1}
  [Q]:H_{Q}^{2}=
  \bigcap_{j=1}^{r_{1}}\bigcap_{k=1}^{r_{2}}[P_{1,j},P_{2,k}]:(I_{1,j}I_{2,k})^{\infty}. 
\end{equation}
\end{lemme}
\begin{preuve} i) a) $\Longrightarrow$ b). --- Pour tout couple
  $(i_{0}, j_{0})\in [1,r_{1}]\times[1,r_{2}]$, d'après
  \cite[th.~3]{Boulier2009}, $\{P_{1, j_{0}}, P_{2, k_{0}}\}$ est un
  ensemble caractéristique de l'idéal $[P_{1, j_{0}}, P_{2,
      k_{0}}]:(I_{1, j_{0}}I_{2, k_{0}})^{\infty}$, puisque les
  $P_{i,k}$ sont de degré $1$ en $\upsilon_{i}$, ce qui implique leurs
  séparants sont égaux à leurs initiaux et que nos
  hypothèses impliquent aussi que $\{P_{1, j_{0}}, P_{2, k_{0}}\}$ est un
  ensemble caractéristique de l'idéal $(P_{1, j_{0}}, P_{2,
    k_{0}}):(I_{1, j_{0}}I_{2, k_{0}})^{\infty}$.  L'idéal $[P_{1,
      j_{0}}, P_{2, k_{0}}]$ contient $[Q]$.

  Donc, si le reste $R$ de la pseudo-réduction de
  $\Spol(Q_{1}, Q_{2})$
  par $Q$ est non nul, il ne dépend d'aucune dérivée stricte de
  $\upsilon_{1}$ ou de $\upsilon_{2}$ et est alors partiellement
  réduit par rapport à $\{P_{1, j_{0}}, P_{2, k_{0}}\}$ pour tout
  $(j_{0},k_{0})\in[1,r_{1}]\times[1,r_{2}]$. Toujours d'après
  \cite[th.~3]{Boulier2009}, il appartient donc à l'idéal algébrique
  $(P_{1, j_{0}}, P_{2, k_{0}}):(I_{1, j_{0}}I_{2, k_{0}})^{\infty}$ pour tout
  $(j_{0},k_{0})\in[1,r_{1}]\times[1,r_{2}]$.

  Soit $\Xi$ l'ensemble des dérivées présentes dans les polynômes
  $P_{i,k}$, pour $i=1,2$ et $1\le k\le r_{i}$ et différentes des
  dérivées de tête $\upsilon_{1}$ et $\upsilon_{2}$. Soit
  $\cK:=\cF(\Xi)$. Les idéaux $(P_{1, j}, P_{2, k})$ pour tous les
  couples $(j,k)\in [1,r_{1}]\times[1,r_{2}]$ définissent $r_{1}\times
  r_{2}$ points distincts $(\eta_{1,j},\eta_{2,k})$ dans l'anneau
  $\cK[\upsilon_{1},\upsilon_{2}]$. Le reste de la pseudo-réduction de
  $\Spol(Q_{1}, Q_{2})$ par $Q$ est de degré en $\upsilon_{1}$
  strictement inférieur à $r_{1}$. S'il est non nul, on peut donc
  trouver $j_{0}$ tel que son évaluation à la valeur $\eta_{1,j_{0}}$
  de $\upsilon_{1}$ définie par $P_{1, j_{0}}$ soit un polynôme non
  trivial de $\cK[\upsilon_{2}]$ de degré strictement inférieur à
  $r_{2}$. Il existe donc $k_{0}$ tel qu'il ne s'annule pas à la
  valeur $\eta_{2,k_{0}}$ de $\upsilon_{2}$ définie par $P_{2,
    k_{0}}$. Ceci implique que $\Spol(Q_{1}, Q_{2})$ n'appartient pas
  à $(P_{1, j_{0}}, P_{2, k_{0}})$: une contradiction. Ce S-polynôme
  est donc bien réduit à $0$ par $Q$.
\smallskip
  
b) $\Longrightarrow$ a). --- Si $T:=\Spol(P_{1, i_{0}},P_{2, j_{0}})$
n'est pas réduit à $0$ par $\{P_{1, i_{0}}, P_{2, j_{0}}\}$, alors son
reste $R$ appartient à $\cK[\upsilon_{1},\upsilon_{2}]$ et donc
également à $(Q):S_{Q}^{\infty}$ par~\cite[th.~3]{Boulier2009}. Il ne
s'annule pas au point $(\eta_{1,i_{0}},\eta_{2,j_{0}})$ défini par
l'idéal $(P_{1, i_{0}}, P_{2, j_{0}})_{\cK[\upsilon]}$. Par ailleurs,
$[P_{1, i_{0}}, P_{2, j_{0}}]=(Q):(\prod_{i\neq i_{0}} P_{1,
  i}\prod_{j\neq j_{0}} P_{2, j})^{\infty}$, donc il existe
$(r_{1},r_{2})\in\N^{2}$ tel que $R(\prod_{i\neq i_{0}} P_{1,
  i}\prod_{j\neq j_{0}} P_{2, j})^{r_{1}}S_{Q}^{r_{2}}\in(Q)$.

L'assertion~b) implique que tous les polynômes de $(Q)$ sont réduits à
$0$ par l'ensemble $Q$, donc il existe $s\in\N$ tel que
$H_{Q}^{s}R(\prod_{i\neq i_{0}} P_{1, i}\prod_{j\neq j_{0}} P_{2,
  j})^{r_{1}}\in(Q)$ et ainsi $[H_{Q}^{s}R(\prod_{i\neq i_{0}} P_{1,
    i}\prod_{j\neq j_{0}} P_{2,
    j})^{r_{1}}](\eta_{1,i_{0}},\eta_{2,j_{0}})$ est nul, ce qui est impossible,
puisque les idéaux premiers définis par les $P_{i,k}$ sont tous
différents, ce qui implique $H_{Q}\neq0\in\cK$ et que par ailleurs
$$R(\eta_{1,i_{0}},\eta_{2,j_{0}})\neq 0\> \hbox{et}\> (\prod_{i\neq i_{0}}
P_{1, i}\prod_{j\neq j_{0}} P_{2,
  j})(\eta_{1,i_{0}},\eta_{2,j_{0}})\neq0.$$

  ii) Comme $Q$ est autoréduit, il suffit de s'assurer que tout
  polynôme appartient à $(Q):H_{Q}^{\infty}$ ssi il est réduit à $0$
  par $Q$. Le lemme de Rosenfeld~\cite[th.~3]{Boulier2009} permet de
  se ramener au cas où un tel polynôme est partiellement réduit. Comme
  les initiaux $I_{Q_{i}}=\prod_{k=1}^{r_{i}}I_{i,k}$ appartiennent à
  $K[\Xi]$, ceci revient à montrer que $Q$ est une base standard de
  $(Q)_{\cK[\upsilon]}:S_{Q}^{\infty}$. C'est bien une base standard
  de $[Q]_{\cK[\upsilon_{1},\upsilon_{2}]}$, en vertu du premier
  critère de Buchberger~\cite{Buchberger1979}, puisque les monômes de
  tête sont étrangers. Par ailleurs,
  $Q_{\cK[\upsilon]}=Q_{\cK[\upsilon]}:S_{Q}^{\infty}$. En effet,
  puisque $\upsilon_{i}$ n'intervient pas dans $Q_{j}$, pour $\{i,
  j\}=\{1,2\}$, l'ensemble des zéro de $Q$ est de la forme
  $(\eta_{1,k_{1}},\eta_{2,k_{2}})$, pour
  $(k_{1},k_{2})\in[1,r_{1}]\times[1,r_{2}]$. Comme les facteurs des
  $P_{i,k}$ contenant $\upsilon_{i}$ sont tous différents par
  hypothèse, $Q_{i}$ et $S_{Q_{i}}$ sont sans facteur commun, pour $i=1,2$
  et donc $(Q_{i})_{\cK[\upsilon]}:S_{i}^{\infty}=(Q_{i})$. On a ainsi
  $(Q)_{\cK[\upsilon]}:S_{Q}^{\infty}=(Q_{1})_{\cK[\upsilon]}:S_{1}^{\infty}
  +(Q_{2})_{\cK[\upsilon]}:S_{2}^{\infty}
    =(Q_{1})_{\cK[\upsilon]}+(Q_{2})_{\cK[\upsilon]}=(Q)_{\cK[\upsilon]}$.

  Par le relèvement du lemme de Lazard~\cite[th.~4]{Boulier2009},
  l'idéal différentiel $[Q]:H_{Q}^{\infty}$ est radical. Ses
  composantes premières sont les idéaux premiers
  $[P_{1,j},P_{2,k}]:(I_{1,j}I_{2,k})^{\infty}$, pour tous les
  $(j,k)\in[1,r_{1}]\times[1,r_{2}]$, puisque les facteurs des
  $P_{i,k}$ contenant leur dérivées de tête sont les seuls facteurs de
  $Q_{i}$ qui ne divisent pas $S_{Q}$. Ceci prouve
  l'équation~\eqref{eq::1}.
\end{preuve}

\begin{exemple}\label{ex::2} Soit $\cF:=\Q(\omega)$ avec
  $\omega^{3}=1$, muni de 
  dérivations $\delta_{i}$, $1=1,2,3$ qui commutent. Soit
  $\cR:=\cF\{x\}$, on définit
  $P_{i,k}:=\delta_{i}x+\omega^{k}\delta_{3}x$, pour $i=1,2$ et $1\le
  k\le 3$. Ainsi, si
  $Q_{1}:=\prod_{k=1}^{3}P_{1,k}=(\delta_{1}x)^{3}+(\delta_{3}x)^{3}$ et
  $Q_{2}:=\prod_{k=1}^{3}P_{1,k}=(\delta_{2}x)^{3}+(\delta_{3}x)^{3}$, alors
  le S-polynôme de $Q_{1}$ et $Q_{2}$ est
  $$
  3(\delta_{2}x)^{2}\delta_{2}Q_{1}-3(\delta_{1}x)^{2}\delta_{1}Q_{2}=
  9(\delta_{3}x)^{2}((\delta_{2}x)^{2}\delta_{2}\delta_{3}x
  -(\delta_{1}x)^{2}\delta_{1}\delta_{3}x),
$$
qui est pseudo-réduit à $0$ par $Q$, en utilisant le
même ordre sur les dérivées qu'à l'exemple~\ref{ex::1}.
\end{exemple}

Dans le cas où des puissances strictes interviennent dans les
produits, nous auront besoin en complément du lemme technique qui suit.

\begin{lemme}\label{lem:carre} Soit $\cP\subset\cR$ un idéal
  premier. On note $\Omega_{\cR/\cF}$ le module des différentielles de
  Kähler et $\rd$ la dérivation canonique de $\cR$ dans $\Omega_{\cR/\cF}$.

i) Pour tout $r\in\N$ et $Q\in\cP^{r}$, on a: $\rd
Q\in\cP^{r-1}\Omega_{\cR/\cF}$.

ii) Soient $\cA=\{A_{1}, \ldots, A_{s}\}$ un ensemble caractéristique
de $\cP$ pour un ordre $\prec$ et $d_{i}\in\N$, pour $1\le i\le s$, des
entiers quelconques.

\noindent a) Pour tout $A\in\cR$, tout $d\in\N$ et tout
$(\tau,\theta)\in\Theta^{2}$, on a $\partial_{\tau\!\! A}\, \theta\!\! A^{d}\in
          [A^{d-1}]$.

\noindent b) Si $Q\in\cQ:=[A_{i}^{d_{i}}|1\le i\le
  s]:H_{\cA}^{\infty}$, on a: $\rd Q=\sum_{i=1}^{s}\sum_{k\in\N}c_{i,
  k}\rd A_{i}^{(k)}$, avec $c_{i, k}\in
          [A_{\iota}^{\bar d_{\iota}}|1\le\iota\le
            s]:H_{\cA}^{\infty}$, où $\bar d_{\iota}:=d_{i}$ si
          $\iota\neq i$ et $\bar d_{\iota}:=d_{i}-1$ si
          $\iota=i$.

\noindent c) Avec les même notations qu'en b), soit $\Xi$ est l'ensemble
des dérivées apparaissant dans les $A_{i}$ et différentes des dérivées
de tête $\upsilon_{A_{i}}$. Si les $A_{i}$ sont tous linéaires en leur
dérivée de tête $\upsilon_{i}$ et si $I_{i}\in\cF[\Xi]$, on a
l'égalité $\cQ\cap\cR_{0}=(A_{i}^{d_{i}}|1\le i\le s):I_{\cA}^{\infty}$,
avec $\cR_{0}:=\cF[\Xi,\upsilon_{A_{1}}, \ldots,
  \upsilon_{A_{s}}]$. Dans ce cas, $\{A_{i}^{d_{i}}|1\le
i\le s\}$ est un ensemble charactéristique de $(A_{i}^{d_{i}}|1\le
i\le s):I_{\cA}^{\infty}$.
\end{lemme}
\begin{preuve}
i) Par définition, tout élément $Q$ de $\cP^{r}$ peut s'exprimer
comme une combinaison linéaire
$$
  Q=\sum_{j=1}^{p}L_{j}m_{j}(B_{1}, \ldots, B_{q})
$$
de $p$ monômes de degré $r$ dépendant d'une famille $B_{i}$ de $q$
éléments de $\cP$. On a alors
$$\begin{array}{ll} \rd Q&=\sum_{j=1}^{p}\left(m_{j}(B_{1}, \ldots, B_{q})\rd L_{j}
  (\in\cP^{r}\Omega_{\cR/\cF})\right.\cr
  &\phantom{\sum_{j=1}^{p}\left(\right.}+\left.L_{j}\sum_{i=1}^{q}\frac{\partial
  m_{i}}{\partial B_{i}}(B_{1}, \ldots,
  B_{q})\rd B_{i}(\in\cP^{r-1}\Omega_{\cR/\cF})\right).
\end{array}  
$$

ii) a) On procède par récurence sur l'ordre de $\theta$. S'il et nul,
c'est-à-dire si $\theta=1$, le résultat est immédiat. Sinon,
$\theta=\delta_{i_{0}}\hat\theta$, avec
$\ord\hat\theta=\ord\theta-1$. Supposons le résultat vrai pour tout
$\vartheta\in\Theta$ d'ordre strictement inférieur à $\ord\theta$.

Pour tout polynôme
$R\in\Q_{\Delta}\{A\}$, on a:
$$
\delta_{i_{0}}R=\sum_{\bar\theta\in\Theta}\delta_{i_{0}}\bar\theta\!\!
A\frac{\partial}{\partial\,\bar\theta\!\! A}R.
$$
Donc, si $\delta_{0}$ ne divise pas $\tau$, on a
$$
\frac{\partial\,\theta\!\! A^{d}}{\partial\,\tau\!\! A}=
\frac{\partial}{\partial\,\tau\!\! A}\sum_{\bar\theta\in\Theta}\delta_{i_{0}}\bar\theta
A\frac{\partial}{\partial\, \bar\theta\!\! A}\hat\theta A^{d}=
\delta_{i_{0}}\frac{\partial\,\hat\theta\!\! A^{d}}{\partial\,\tau\!\! A}.
$$
Sinon,
on a $\tau=\delta_{i_{0}}\hat\tau$ et on peut écrire:  
$$
\frac{\partial\,\theta A^{d}}{\partial\,\tau\!\! A}=
\frac{\partial}{\partial\,\tau\!\! A}
\sum_{\bar\theta\in\Theta}\delta_{i_{0}}\bar\theta\!\!
A\frac{\partial}{\partial\, \bar\theta\!\! A}\hat\theta A^{d}=
\delta_{i_{0}}\frac{\partial\,\hat\theta\!\! A^{d}}{\partial\,\tau\!\! A}+
\frac{\partial\,\hat\theta\!\! A^{d}}{\partial\,\hat\tau\!\! A},
$$ où $\partial\,\hat\theta\!\! A/\partial\,\tau\!\! A$ et $\partial\,\hat\theta\!\!
A/\partial\,\hat\tau\!\! A$ appartiennent à $[A^{d-1}]$, en utilisant
l'hypothèse de récurrence, ce qui implique le résultat.

b) Si $Q\in\cQ$, alors il existe $h\in\N$ tel que
$$
I_{\cA}^{h}Q=\sum_{i=1}^{s}\sum_{\theta\in\Theta}
L_{i,\theta}\theta\!\!A_{i}^{d_{i}},
$$
donc $\rd\left(I_{\cA}^{h}\right)Q(\in\cQ\Omega_{\cR/\cF})+I_{\cA}^{h}\rd
Q=$

$$\begin{array}{l}
  \sum_{i=1}^{s}\sum_{\theta\in\Theta}
\theta\!\!A_{\iota}^{d_{\iota}}\rd L_{\iota,\kappa}
(\in\cQ\Omega_{\cR/\cF})\\
+ \sum_{i=1}^{s}\sum_{\theta\in\Theta} \sum_{\tau\in\Theta}
L_{i,\theta}\frac{\partial\,
  \theta\!\!A_{i}^{d_{i}}} {\partial\, \tau\!\! A_{i}}\rd
\tau A_{i} (\in [A_{i}^{d_{i}-1}]\Omega_{\cR/\cF}\> \hbox{par a),}
\end{array}
$$
d'où le résultat.

c) L'inclusion $\supset$ est immédiate. Pour $\subset$, on procède par
récurrence sur $d\in(\N^{\ast})^{s}$, en utilisant l'ordre sur
$(\N^{\ast})^{{s}}$ défini par $\hat d\le \tilde d$ si $\hat d_{i}\le
\tilde d_{i}$ pour $1\le i\le s$. Pour $(d_{1}, \ldots, d_{s})=(1,
\ldots, 1)$, le résultat est vrai, puisque $\cA$ est un ensemble
caractéristique de $\cA$. Soit $Q\in\cQ\cap\cR_{0}$, nous allons
prouver que $Q\in(A_{i}^{d_{i}}|1\le i\le s)_{\cR_{0}}$, en supposant
le résultat vrai pour tout $\hat d<d$.

D'après le premier critère de Buchberger~\cite{Buchberger1979},
l'ensemble $\{A_{i}^{d_{1}}, \ldots, A_{s}^{d_{s}}\}$ est une base
standard algébrique de l'idéal qu'il engendre, dans l'algèbre
$\cF(\Xi)[\upsilon_{A_{1}}, \ldots, \upsilon_{A_{s}}]$. C'est donc
aussi un ensemble caractéristique de l'idéal qu'il engendre dans
$\cF[\Xi,\upsilon_{A_{1}}, \ldots, \upsilon_{A_{s}}]$.

Si $Q\notin (A_{i}^{d_{i}}|1\le i\le
s)_{\cR_{0}}$, qu'il est non nul et irréductible pour la réduction par
l'ensemble caractéristique algébrique $\{A_{1}^{d_{1}}, \ldots,
A_{s}^{d_{s}}\}$, le degré de $Q$
en $\upsilon_{A_{i}}$ est strictement inférieur à $d_{i}$.

Comme $Q\in(\cA):H_{\cA}^{\infty}$, $Q$ est réductible par $\cA$ et dépend donc
d'une des dérivées de tête $\upsilon_{A_{i_{0}}}$; nécessairement, on
a $d_{i_{0}}>1$, sinon $Q$ serait réductible par
$A_{i_{0}}^{d_{i_{0}}}$. Par b), le coefficient de $\rd A_{i_{0}}$
dans $\rd Q$ appartient à $[A_{i}^{\hat d_{i}}]:H_{\cA}^{\infty}$,
avec $\hat d_{i}=d_{i}$ si $i\neq i_{0}$ et $\hat
d_{i_{0}}=d_{i_{0}}-1$. Comme le degré de $Q$ en $\upsilon_{i}$ est au
plus $d_{i}-1$, ce coefficient qui est dans $\cR_{0}$ est de degré en
$\upsilon_{i}$ au plus $\hat d_{i_{0}}-1$ et irréductible par
$\{A_{i}^{\hat d_{i}}|1\le i\le s\}$ ce qui contredit, s'il est non
nul, l'hypothèse de récurrence.

Tout élément de $\cQ$ est donc réduit à $0$ par $\{A_{i}^{\hat
  d_{i}}|1\le i\le s\}$, qui est autoréduit et constitue donc un ensemble
charactéristique de $\cQ$.
\end{preuve}

\begin{exemple}\label{ex::3} On prend $\cF:=\Q$, muni de $3$
  dérivations $\delta_{i}$, $i=1,2,3$ et $\cR:=\cF\{x\}$, avec le même
  ordre sur les dérivées qu'à l'exemple~\ref{ex::1}. On prend
  $\cP:=[\delta_{1}x,\delta_{2}x,\delta_{3}x]$. Alors, $\rd(
  \delta_{1}x\delta_{2}x\delta_{3}x(\in\cP^{3}))=\delta_{2}x\delta_{3}x\delta_{1}\rd
  x + \delta_{1}x\delta_{3}x\delta_{2}\rd
  x  +\delta_{1}x\delta_{2}x\delta_{3}\rd
  x\in\cP^{2}\cM$. De même,
  $\rd
  [(\delta_{1}x)^{2}+(\delta_{2}x)^{3}+
    (\delta_{3}x)^{4}]=$
  $$2\delta_{1}x\delta_{1}\rd x (\in [\delta_{1} x]\cM)
  +3(\delta_{2}x)^{2}\delta_{2}\rd x (\in [\delta_{2} x]^{2}\cM)+
  4(\delta_{3}x)^{3}\delta_{3}\rd x(\in [\delta_{3} x]^{3}\cM)
  $$
  et
  $\cB:=\{(\delta_{1}x)^{2},(\delta_{2}x)^{3},(\delta_{3}x)^{4}\}$ est un
  ensemble caractéristique de l'idéal qu'il engendre, mais pas une
  représentation caractéristique, puisque $\delta_{1}^{2}x$ est réduit
  à $0$ par $\cB$ sans appartenir à $\cB$. 
\end{exemple}
\medskip

On peut maintenant énoncer cette nouvelle version du théorème de
Hashemi et Touraji~\cite{Hashemi2014}. La différence essentielle par
rapport à la version originale, outre l'utilisation d'une version un
peu étendue du critère de Boulier, est que l'on concentre tous les
facteurs de $Q_{i}$ dont la dérivée dominante est strictement
inférieure à la dérivée dominante commune aux facteurs linéaires
homogènes en $z$, en un facteur unique $Q_{i, 0}$, sur lequel il n'est
pas fait d'autres hypothèses. On notera que c'est un facteur du
séparant $S_{Q_{i}}$, qui est implicitement supposé non nul dans i), ce qui
exclut les zéros de ce terme.

\begin{theoreme}\label{th:principal} On pose
  $\bar Q_{i}=\prod_{k=0}^{r_{i}} P_{i,  k}^{d_{i,  k}}$ et $Q_{i}:=\prod_{k=1}^{r_{i}}
  P_{i, k}^{d_{i, k}}$, $i=1, 2$, où pour tout 
  couple $(P_{1, k_{1}},P_{2,  k_{2}})$, $(k_{1}, k_{2})\in[1,
    r_{1}]\times[1, r_{2}]$,
  les hypothèses du th.~\ref{th:Boulier} sont satisfaites pour un même
  $\theta_{i}z$, et $\upsilon_{P_{i, 0}}\prec \theta_{i}\upsilon_{z}$,
  $i=1, 2$. On suppose en outre que les $P_{i,  k}$ sont tous différents.

  i) Sous ces hypothèses, et si $d_{i,  k}=1$, pour $i=1, 2$ et $1\le k\le
  r_{i}$, alors
$\bar Q:=\{\bar Q_{1}, \bar Q_{2}\}$ est un ensemble caractéristique et une
  représentation caractéristique de
$[\bar Q]:S_{\bar Q}^{\infty}=[Q]:S_{Q}^{\infty}$.

ii) Sous ces hypothèses, et pour des $d_{i,  k}$ quelconques, 
$
T:=\Spol(\bar Q_{1},\bar Q_{2})= S_{\bar Q_{2}}\theta_{2}\bar Q_{1} -
S_{\bar Q_{1}}\theta_{1}\bar Q_{2},
$
est pseudo-réduit à $0$ par $\bar Q$.
\end{theoreme}
\begin{preuve} i) D'après le th.~\ref{th:Boulier},
  tous les S-polynômes pour tous les couples sont réduits à $0$ et
  d'après le lem.~\ref{lem:prod}, $\{Q_{1}, Q_{2}\}$ est
  un ensemble caractéristique de $[Q_{1},
    Q_{2}]:S_{Q}^{\infty}$. Donc, les $\bar Q_{i}$
  qui sont des multiples et ont les mêmes dérivées de têtes
  respectives que les $Q_{i}$, forment un ensemble
  caractéristique du même idéal (et réduisent donc à $0$ leur
  S-polynôme). On a alors $S_{\bar Q}=P_{1,0}P_{2,0}S_{Q}$ et
  $I_{\bar Q}=P_{1,0}P_{2,0}I_{Q}$, donc
  $[\bar Q]:H_{\bar Q}^{\infty}=[Q]:H_{Q}^{\infty}=[Q]:S_{Q}^{\infty}$.

  ii) Si $T\in[Q_{1},Q_{2}]$ contient une dérivée stricte de
  $\theta_{i}\upsilon_{i}$, pour $i=1$ ou $i=2$, alors $T$ est
  réductible par $\{Q_{1}, Q_{2}\}$. Sinon, pour tout couple $(k_{1},
  k_{2})\in[1, r_{1}]\times[1, r_{2}]$, le
  lemme~\ref{lem:carre}~ii)~c) implique $T\in(P_{1,
    k_{1}}^{d_{1,k_{1}}},P_{2, k_{2}}^{d_{2,k_{2}}}):S_{z}^{\infty}$,
  puisque $I_{P_{i,k}}$ est égal à $S_{z}$. On a donc
$$
T\in\bigcap_{(k_{1},k_{2})\in[1, r_{1}]\times[1, r_{2}] }\left(P_{1, k_{1}}^{d_{1,k_{1}}},P_{2,
  k_{2}}^{d_{2,k_{2}}}\right):S_{z}^{\infty}.
$$ Montrons que $T$ est pseudo-réduit à $0$ par $\{Q_{1},Q_{2}\}$, ce
qui revient à montrer qu'il est réduit à $0$, au sens des bases
standard, par les $Q_{i}$, considérés comme des polynômes de
$\cF(\Xi)[\upsilon_{1},\upsilon_{2}]$, où $\Xi$ est défini comme au
lem.~\ref{lem:carre}~ii)~c). Dans cet anneau, les éléments $R$ de
l'idéal primaire $(P_{1, k_{1}}^{d_{1,k_{1}}},P_{2,
  k_{2}}^{d_{2,k_{2}}})$ sont tels que, pour tout couple d'entiers
$(\alpha_{1},\alpha_{2})$, tels que $0\le \alpha_{i}<d_{i,k_{i}}$,
$\partial_{\upsilon_{1}}^{\alpha_{1}}\partial_{\upsilon_{2}}^{\alpha_{2}}R$
appartienne à l'idéal premier $(P_{1, k_{1}},P_{2,
  k_{2}}):S_{z}^{\infty}$, qui définit le point
$(\eta_{1,k_{1}},\eta_{2,k_{2}})$\footnote{C'est un cas particulier
  simple de la description d'un idéal algébrique primaire par un
  système d'opérateurs différentiels linéaires. Voir Cid-Ruiz \textit{et
    al.}~\cite{Cid-Ruiz2021} pour plus de détails et de
  références.}. Si
$R\in[Q_{1},Q_{2}]$ est le reste de la réduction de $T$ et s'il est
non nul, c'est un polynôme de degré $\beta_{i}$ en $\upsilon_{i}$
strictement inférieur à $\sum_{i=1}^{r_{i}}d_{i,k_{i}}$, pour $i=1$ et
$i=2$. On peut l'écrire comme un polynôme en $\upsilon_{1}$ avec des
coefficients polynômes en $\upsilon_{2}$, c'est-à-dire:
$R=\sum_{i=0}^{\beta_{1}}c_{i}(\upsilon_{2})\upsilon_{1}^{i}$. Comme
$R$ est non nul, le coefficient $c_{\beta_{1}}$ est un polynôme non
nul de degré strictement inférieur à $\sum_{i=1}^{r_{2}}d_{2,k_{2}}$,
qui ne peut donc s'annuler en les $r_{2}$ points $\eta_{2,k_{2}}$,
pour $1\le k_{2}\le r_{2}$, avec mutiplicité $d_{2,k_{2}}$. Soit
$k_{0}$ tel que la mutiplicité en $\eta_{2,k_{0}}$ soit
$\gamma<d_{2,k_{0}}$. Alors
$\partial_{\upsilon_{2}}^{\gamma}c_{\beta}$ est un élément non nul du
corps $\cF(\Xi)$. On doit avoir $(\partial_{\upsilon_{1}}^{\alpha_{1}}
\partial_{\upsilon_{2}}^{\gamma}R)(\eta_{1,k_{1}},\eta_{2,k_{0}})=0$,
pour tout $1\le k_{1}\le r_{1}$ et tout $0\le\alpha_{1}<d_{1,k_{1}}$,
ce qui est impossible, puisque l'ordre de $R$ en $\upsilon_{1}$ est
strictement inférieur à $\sum_{i=1}^{r_{1}}d_{1,k_{1}}$. On en déduit
que $R$ est nul, ce qui achève la preuve, puisque réduire par $\bar Q$
ou par $Q$ est équivalent, à une multiplication près par une puissance
de $P_{1,0}P_{2,0}$.
\end{preuve}
\medskip

\begin{exemple}\label{ex::4}
On prend le corps $\cF:=\Q$, équipé de trois dérivations $\delta_{1}$,
$\delta_{2}$, $\delta_{3}$
et
$\cR:=\cF\{x\}$ avec le même ordre sur les dérivées qu'à
l'exemple~\ref{ex::1}.

Soient $Q_{i}:=\prod_{k=1}^{r_{i}}(\delta_{1}x - a_{i,k}\delta_{3}x)$,
pour $i=1,2$, avec $a_{i,k}\in\Q$, $1\le k\le r_{i}$. Alors,
$\{Q_{1},Q_{2}\}$ réduit à $0$ le S-polynôme de $Q_{1}$ et $Q_{2}$ et
est un ensemble caractéristique de l'idéal qu'il engendre. Ceci
demeure vrai
dans les cas où apparaissent des facteurs multiples, c'est-à-dire
quand certains coefficients $a_{i,k}$ sont égaux entre eux.

On peut multiplier $Q_{i}$ par n'importe quel facteur $P_{i,0}$, dont
la dérivée de tête est inférieure à $\delta_{1}x$ et $\delta_{2}x$, par exemple
$\delta_{3}x$. Un tel facteur ne fait que contribuer aux initiaux et
aux séparant lors des réductions.
\end{exemple}

\bigskip

\section*{Conclusion}

Du point de vue de la complexité algorithmique, le nombre de monômes
obtenus en développant les produits rend génériquement la taille de
leur représentation dense supérieure à celle de tous les idéaux
premiers obtenus en factorisant et la version avec carrés n'est pas
utile pour l'algorithme Rosenfeld -- Gröbner~\cite{Boulier2009}.

Ces constats ne doivent pas occulter l'intérêt théorique et pratique
qu'il peut y avoir à considérer des produits et des puissances
d'idéaux en algèbre différentielle et donc à renforcer l'outillage
théorique disponible pour leur étude. On peut par exemple mettre en
œuvre des formes de remontées de Newton--Hensel, ce qui a déjà été mis
à profit pour calculer des solutions en séries~\cite{Bostan2007}. La
méthode pourrait s'étendre pour permettre le calcul d'un ensemble
caractéristique pour un système dont une solution régulière est
connue, ou aussi une représentation des solutions utilisant la
représentation compacte des polynômes différentiels comme des
programmes d'évaluation, dans l'esprit de l'algorithme de
Kronecker~\cite{Giusti2001,DAlfonso2011}.

\end{document}